\documentclass{article}
\usepackage{multicol}

\usepackage{amsmath,amsfonts}
\usepackage{graphicx}
\usepackage{tikz}
\usepackage{verbatim}
\usepackage{color}
\usepackage{comment}
% ...

\newtheorem{thm}{Theorem}

\newtheorem{fct}[thm]{Fact}
\newtheorem{lem}[thm]{Lemma}

\def\s{{\sigma}}

\def\zN{{\mathbb N}}

\def\cF{{\cal F}}
\def\cG{{\cal G}}
\def\cH{{\cal H}}

\def\Pf{{\noindent\bf Proof.\ \ }}

\def\pf{{\hfill $\Box$\bigskip}}
\def\ra{{\rightarrow}}

\def\pot{{\sf pot}}

\definecolor{brwn}{RGB}{140, 70, 20}
\definecolor{gren}{RGB}{  0,140, 10}

\newcommand{\la}[1]{\textcolor{blue}{\sf{#1}}}
\newcommand{\mg}[1]{\textcolor{gren}{\sf{#1}}}
\newcommand{\gh}[1]{\textcolor{brwn}{\sf{#1}}}

\begin{document}

\title{$t$-Pebbling in $k$-connected diameter two graphs} 

\author{
Liliana Alc\' on
\and
Marisa Gutierrez
\and
Glenn Hurlbert 
}
\date{}
\maketitle

\begin{abstract}
Graph pebbling models the transportation of consumable resources.
As two pebbles move across an edge, one reaches its destination while the other is consumed.
The $t$-pebbling number is the smallest integer $m$ so that any initially distributed supply of $m$ pebbles can place $t$ pebbles on any target vertex via pebbling moves.
The 1-pebbling number of diameter two graphs is well-studied.
Here we investigate the $t$-pebbling number of diameter two graphs under the lense of connectivity.
\end{abstract}

\section{Introduction}\label{s:intro}

Graph pebbling has an interesting history, with many challenging open problems.
Calculating pebbling numbers of graphs is a well known computationally difficult problem.
See \cite{Hurl1,Hurl2} for more background.

A {\it configuration} C of pebbles on the vertices of a connected graph G is a function $C : V(G)\ra \zN$ (the nonnegative integers), so that $C(v)$ counts the number of pebbles placed on the vertex $v$. 
We write $|C|$ for the {\it size} $\sum_v C(v)$ of $C$; i.e. the number of pebbles in the configuration. 
A {\it pebbling step} from a vertex $u$ to one of its neighbors $v$ reduces $C(u)$ by two and increases $C(v)$ by one. 
Given a specified {\it root} vertex $r$ we say that $C$ is $t$-{\it fold} $r$-{\it solvable} if some sequence of pebbling steps places $t$ pebbles on $r$. 
We are concerned with determining $\pi_t(G,r)$, the minimum positive integer  $m$ such that every configuration of size $m$ on the vertices of $G$ is $t$-fold $r$-solvable. 
The $t$-{\it  pebbling number} of $G$ is defined to be $\pi_t(G) = \max_{r\in  V(G)} \pi(G,r)$.
We avoid  writing $t$ when  $t=1$.

Pebbling number of diameter 2 graphs was solved and characterized
by the following theorem.
For the purpose of the present work, it is enough to know that a pyramidal graph has no {\it universal} vertex (a vertex adjacent to every other vertex) and has connectivity 2. 
\begin{thm}
\label{t:Diam2}
\cite{ClHoHu,PaSnVo}
For a diameter 2 graph $G$  with connectivity $k$ and $n$ vertices,  $\pi(G)=n+1$ if and only if $k=1$ or
$G$ is pyramidal. Otherwise (i.e. $k=2$ and $G$ is not pyramidal, or $k\geq3$), $\pi(G)=n$.
\end{thm}

In contrast, other than the following bound,
little is known about the $t$-pebbling number of diameter 2 graphs. 

\begin{thm}\label{t:d2bound}\cite{HeHeHu}
If $G$ is a diameter 2 graph on $n$ vertices then $\pi_t(G) \le \pi(G)+4t-4$. Moreover, $\liminf_{t\ra\infty}\pi_t(G)/t = 4$.
\end{thm}
        
The goal of the present paper is to determine the exact $t$-pebbling number of a large 
 subfamily of diameter 2 graphs by considering their connectivity.
Define $\cG(n,k)$ to be the set of all $k$-connected graphs on $n$ vertices having a universal vertex.
Set $f_t(n,k) = n + 4t - k - 2$ and $h_t(n)=n+2t-2$.
Notice that $h_t(n)\ge f_t(n,k)$ if and only if $k\ge 2t$.  
Define $p_t(n,k)=\max\{f_t(n,k),h_t(n)\}$.
 The main result is the following theorem which is proved in Section \ref{s:thm}. 
\begin{thm}\label{Thm}
If $G\in\cG(n,k)$ then $\pi_t(G)=p_t(n,k)$.
\end{thm}

We observe from our result that, for any fixed $t$,
in the family of graphs with universal vertex, there are graphs whose $t$-pebbling number is much lower than
the bound given by Theorem  \ref{t:d2bound}, and also that there are graphs reaching that bound:
when $k\ge 2t$ we have $\pi_t(n,k)=(n+4t-4)-2(t-1)$; when $k< 2t$ $\pi_t(n,k)=(n+4t-4)-(k-2)$.


It will be useful to take advantage of Menger's Theorem.
The version of Menger's theorem that we use is the following (exercise 4.2.28 in \cite{West}).

\begin{thm}
\label{Menger}
{\bf (Menger's Theorem)}
\cite{West}
Let $G$ be a $k$-connected graph and $S=\{v_1,\ldots,v_k\}$ be a multiset of vertices of $G$.
For any $r\not\in S$ there are $k$ pairwise-internally-disjoint paths, one from each $v_i$ to $r$.
\end{thm}

\section{Technical Lemmas} \label{s:tech}


We begin with a lemma that is used to prove lower bounds on the pebbling number of a graph by helping to show that certain configurations are unsolvable.

For a vertex $v$, define its {\it open neighborhood} $N(v)$ to be the set of vertices adjacent to $v$, and its {\it closed neighborhood} $N[v]=N(v)\cup\{v\}$.
We say that a vertex $y$ is a {\it junior sibling of} a vertex $x$ (or, more simply, {\it junior to} $x$) if $N(y)\subseteq N[x]$, and that $y$ is a {\it junior} if it is junior to some vertex $x$.

\begin{lem}
\label{l:JuniorRemoval}
{\bf (Junior Removal Lemma)}
\cite{AGH3}
Given the graph $G$ with root $r$ and $t$-fold $r$-solvable configuration $C$, suppose that $y$ is a junior with $C(y)=0$.
Then $C$ (restricted to $G-y$) is $t$-fold $r$-solvable in $G-y$.
\end{lem}


Given a configuration $C$  of pebbles, we say that a path $Q=(r,q_1,\ldots,q_j)$  with $j\geq 1$ is a {\it slide} from $q_j$ to $r$ if no $q_i$ is zero (it has no pebbles on)  and $q_j$ has at least two pebbles.

A {\it potential move} is a pair of pebbles sitting on the same vertex.
To say that $C$ has $j$ potential moves means that the $j$ pairs are pairwise disjoint.
For example, any configuration on 5 vertices with values $0,1,1,2,$ and $7$ has 4 potential moves.
The {\it potential} of $C$, $\pot(C)$, is the maximum $j$ for which $C$ has $j$ potential moves.
Because every solution that requires a pebbling move uses a potential move, the following fact is evident. 

\begin{fct} \label{f:pot}
Let $r$ be an empty vertex in a configuration $C$ with $\pot(C)<t$. 
Then $C$ is not $t$-fold $r$-solvable.
\end{fct}

Basic counting yields the following lemma.

\begin{lem}\label{l:PotLem}
{\bf (Potential Lemma)}
Let $G$ be a graph on $n$ vertices. 
If $C$ is a configuration on $G$ of  size $n+y$ ($y\ge 0$) having $z$ zeros, then $\pot(C)\ge\lceil \frac{y+z}{2}\rceil$.
\end{lem}

A nice application of the Potential Lemma is the following, which we will use repeatedly in the arguments that follow.

\begin{lem}\label{l:zeros}
{\bf (Slide Lemma)}
Let $r$ be a vertex of a $k$-connected graph $G$.
Let $C$ be a configuration on $G$ of size $n+y$ ($y\geq 0$) with $z$ zeros.
If $\lceil\frac{y+3z}{2}\rceil\le  k$ then $C$ is $\lceil\frac{y+z}{2}\rceil$-fold $r$-solvable.
\end{lem}

\Pf
Set $p=\lceil\frac{y+z}{2}\rceil$.
By Lemma \ref{l:PotLem} we can choose a set $P$ of $p$ potential moves.
Note that the hypothesis implies that $p\le k-z$.
Delete all non-root zeros to obtain $G'$.
Since $G$ is $k$-connected, $G'$ is $p$-connected.
Thus Menger's Theorem \ref{Menger} implies that there are $p$ pair-wise disjoint slides in $G'$ from $P$ to $r$, which yield $p$ $r$-solutions.
\pf

\section{Proof of Theorem \ref{Thm}} \label{s:thm}

The proof will follow from Lemmas \ref{l:LowerBd} and \ref{l:upper_bound}, below. 
Let $u$ be a universal vertex of a graph $G\in\cG(n,k)$.
If $C$ is a configuration of size $n+2t-3$ with $u$ empty and every other vertex odd then $\pot(C)=t-1$,
and so $C$ is not $t$-fold $u$-solvable.
Hence $\pi_t(G,u)\ge n+2t-2$.
On the other hand, if $|C|\ge n+2t-2$ then $\pot(C)\ge t$ when $u$ is empty, and $\pot(C)\ge t-1$ when $u$ is not; either way $C$ is $t$-fold $u$-solvable because $u$ is universal.
Thus $\pi_t(G,u)= n+2t-2$, which is at most $p_t(n,k)$ always.


\subsection{Lower bound}\label{ss:lower_bound}


Clearly, $\pi_t(G)\geq \pi_t(G,u)=h_t(n)$. Now let $r$ be any non-universal vertex of $G$, and let $s$ be a vertex at distance two from $r$.
Let $X$ be any $(r,s)$-cutset of size $k$ (in particular, $u\in X$) and define the configuration $F_t(n,k)$ by placing 0 on $r$ and $X$, $4t-1$ on $s$, and 1 on each vertex of $V(G)-(X\cup\{r,s\})$; then $|F_t(n,k)| =(4t-1)+(n-k-2)= f_t(n,k)-1$.

Since the vertices of $X-\{u\}$ have 0 pebbles and all them are juniors to $u$, Lemma \ref{l:JuniorRemoval} states that if $t$ pebbles can reach $r$ then $2t$ pebbles can reach $u$. 
But, with exactly $2t-1$ potential moves in $F$, by Fact \ref{f:pot}, we can place at most $2t-1$ pebbles on $u$. Therefore $\pi_t(G,r)\geq f_t(n,k) $, implying  $\pi_t(G)\geq f_t(n,k) $.

We record these results as

\begin{lem}\label{l:LowerBd}
For $G\in\cG(n,k)$ we have $\pi_t(G)\ge p_t(n,k)$.
\end{lem}

\subsection{Upper bound}\label{ss:upper_bound}

We will prove that any configuration of size $f_t(n,k)$ when $k\leq 2t$, and of size $h_t(n)$ when $k\geq 2t$, is $t$-fold $r$-solvable for any $r\in V(G)$.

\begin{lem}\label{l:upper_bound}
For $k\ge 2$, let $G\in\cG(n,k)$ with a universal vertex $u$, and let $r$ be any root vertex.
Then $\pi_t(G,r) \leq p_t(n,k).$
\end{lem}

\begin{figure}
\begin{center}
\includegraphics[height=2.05in]{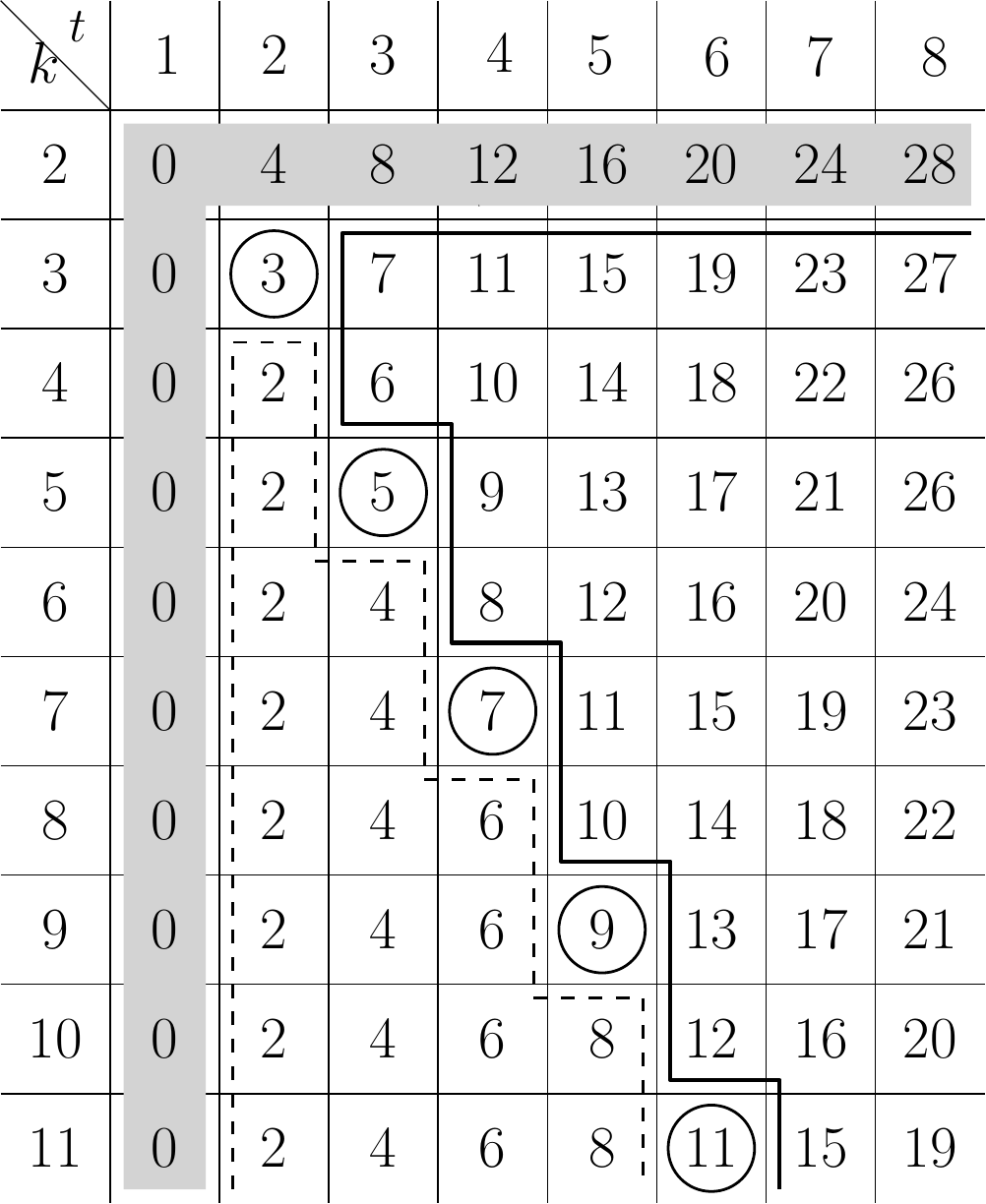}
\end{center}
\caption{The values $m$ for which 
$\pi_t(G)=|V(G)|+m$.} \label{f:chart2}
\end{figure}

\Pf
First note that the lemma is true when $t=1$.
Indeed, in this case we have $k\ge 2t$, and so $p_t(n,k)=h_t(n)=n+2t-2=n$.
On the other hand, because no pyramidal graph has a universal vertex, we have from Theorem \ref{t:Diam2} that $\pi(G)=n$, hence $\pi(G,r)\leq n$.

In addition, the lemma holds for $k=2$. 
Indeed, in this case we have $k\le 2t$, and so $p_t(n,k)=f_t(n,k)=n+4t-k-2=n-4t-4$.
Also, we have by Theorem \ref{t:d2bound} that $\pi_t(G,r)\le n+4t-4$.

Hence, we may assume that $t\ge 2$ and $k\geq 3$.
Figure \ref{f:chart2} shows the structure of this proof.
As was noted above, the 
grey section has been proven before.
We continue by proving the 
dashed-bordered, lower left section and diagonal circled entries together, and then the 
solid-bordered, upper right
section by induction.
\\

\noindent
{\it Base case.}

We will simultaneously address the case $k=2t-1$ (the 
circled entries), for which $|C|=f_t(n,k)=n+2t-1$, and the case  $ k \geq 2t$ (the 
dashed-bordered section), for which
$|C|=h_t(n)=n+2t-2$, by  writing $ k \geq 2t-1  $ and considering a configuration
of size $|C|=n+2t-2+\phi$, where $\phi=1$ if $2t-1=k$ and 0 otherwise.
The natural idea we leverage here is repeating the argument that zeros force potential which, combined with connectivity, yields either more solutions or more zeros.

Let $x \geq 0$ such that $k=2t-1+x$. By Lemma \ref{l:PotLem}, since we may assume that $C(r)=0$ (otherwise induct on $t$), we have at least $\lceil (2t-2+1)/2\rceil=t$ potential moves.
Therefore, we have at least $t$ solutions if there are at least $t$ different slides from them to $r$. 

Thus we consider the case in which there are at most $t-1$ slides; that is, from some of the vertices in which a potential move is sitting, say $v$, there is no path to $r$ without an internal zero after considering the remaining $t-1$ slides. 
Since $G$ is $k$-connected, that implies that $C$ has at least $k-(t-1)$ zeros between $v$ and $r$ and so, because of $r$, $C$ has at least $k-(t-1)+1=t+1+x$ zeros.

Assume that there are exactly $z=t+1+j$ zeros, for some $j\ge x$.
Then, by Lemma \ref{l:PotLem}, $C$ has at least
$$\left\lceil\frac{(2t-2)+(t+1+j)}{2}\right\rceil
=t+\left\lceil\frac{t-1+j}{2}\right\rceil$$
potential moves.
If there are at least $t-\left\lceil\frac{t-1+j}{2}\right\rceil$ slides from them to $r$, then we can use those slides for that many solutions. 
Then, the other $\left\lceil\frac{t-1+j}{2}\right\rceil$ solutions  can be obtained from the remaining $2\left\lceil\frac{t-1+j}{2}\right\rceil$ potential moves, putting $2\left\lceil\frac{t-1+j}{2}\right\rceil$ pebbles on the universal vertex $u$ and then $\left\lceil\frac{t-1+j}{2}\right\rceil$ on $r$.

Otherwise, there are at most $t-\left\lceil\frac{t-1+j}{2}\right\rceil-1$ slides, from which we find, using $k=2t-1+x$, at least
$$k-\left(t-\left\lceil\frac{t-1+j}{2}\right\rceil-1\right)+1
=t+x+\left\lceil\frac{t-1+j}{2}\right\rceil+1$$
zeros.
Clearly, this number cannot exceed the total number of zeros $z=t+1+j$; therefore
$j\ge x+\left\lceil\frac{t-1+j}{2}\right\rceil \ge x+\frac{t-1+j}{2}$, and so $j\ge t-1+2x$.

Let $j=t-1+2x+i$ for some $i\ge 0$; then $z=t+1+j=t+1+t-1+2x+i= 2t+2x+i$. 
Applying Lemma \ref{l:PotLem} again, there are at least
$$\left\lceil\frac{(2t-2)+(2t+2x+i)}{2}\right\rceil=2t+x-1+\lceil i/2\rceil$$ 
potential moves.

If either $x\geq 1$ or $i\ge 1$, then we can move $2t$ pebbles to the universal vertex $u$, and then $t$ to $r$.

Hence, we consider the case for which $x=i=0$; i.e. $k=2t-1$, $z=2t$, and $|C|=n+2t-1$ (because $\phi=1$ in such a case).
We let $T$ be the star centered on $u$, having leaves $r$ and the nonzero vertices of $G$. 
Clearly, $T$ is a subgraph of $G$ with $n+2t-1$ pebbles on it and  with either $2+(n-z)$ or $1+(n-z)$ vertices, depending on whether $u$ is empty or not. 
In either case $n(T)\leq 2+n-z=2+n-2t$. 
Therefore, since
$$    \pi_t(T,r)
    = n(T)+4t-3
    \leq (2+n-2t)+4t-3
       = n+2t-1
    = |C(T)|,$$
we see that $C$ is $r$-solvable.
\\

\noindent
{\it Induction step.}

Finally, we consider the cases when $k<2t-1$ (the 
solid-bordered
section); so $|C|=f_t(n,k)=n + 4t - k - 2$. 
Since $2(t-1)=2t-1-1\geq k$, we have $\pi_{t-1}(G,r)= f_{t-1}(n,k)=n+4(t-1)-k-2=n+4t-k-2-4=|C|-4$. 
Hence, if $C$ has a solution of cost at most 4, we are done.
Otherwise, there is at most one vertex $v$ having two or more pebbles, and on such a vertex there are at most 3 pebbles. 
This implies the contradiction $|C|\leq 3+(n-2)$, which completes the proof.
\pf



In future work we intend to study $k$-connected diameter two graphs without a universal vertex, and use that work as a base step toward studying graphs of larger diameter.



\begin{tabular}{p{0.5\textwidth}p{0.5\textwidth}}
{\parindent=0pt
\obeylines
Liliana Alc\' on
Centro de Matem\' atica La Plata
Universidad Nacional de La Plata
CONICET, Argentina
liliana@mate.unlp.edu.ar
\par
}
&
{\parindent=0pt
\obeylines
Marisa Gutierrez
Centro de Matem\' atica La Plata
Universidad Nacional de La Plata
CONICET, Argentina
marisa@mate.unlp.edu.ar
\par
}\\[1em]
{\parindent=0pt
\obeylines
Glenn Hurlbert
Department of Mathematics and Applied Mathematics
Virginia Commonwealth University
Richmond, Virginia, USA
ghurlbert@vcu.edu
\par
}
&
{\parindent=0pt
\obeylines
\par
}
\end{tabular}

\end{document}